\documentclass[12pt,oneside]{article}

\usepackage{amsfonts}
\usepackage{amscd}
\usepackage{amsthm}
\usepackage{amssymb}
\usepackage{amsmath}
\usepackage{epsfig}
\usepackage{graphicx}

\title{Weakly Hyperbolic Actions of Kazhdan Groups on Tori}

\author{Benjamin Schmidt\footnote{Partially funded by VIGRE grant DMS-9977371}}
\date{}

\theoremstyle{proposition}
\newtheorem{Lem}{Lemma}[section]
\newtheorem{Prop}[Lem]{Proposition}

\newtheorem*{Def}{Definition}

\theoremstyle{plain}

\newtheorem{Thm}[Lem]{Theorem}
\newtheorem{Cor}[Lem]{Corollary}

\theoremstyle{remark}

\newtheorem*{Prf}{Proof}

\DeclareMathOperator{\GL}{GL}
\DeclareMathOperator{\Aut}{Aut}
\DeclareMathOperator{\spn}{span}

\DeclareMathOperator{\Diff}{Diff}
\DeclareMathOperator{\proj}{proj}
\DeclareMathOperator{\length}{length}
\DeclareMathOperator{\Prob}{Prob}
\DeclareMathOperator{\Out}{Out}

\begin{document}
\maketitle

\begin{abstract}
We study the ergodic and rigidity properties of weakly hyperbolic
actions. First, we establish ergodicity for $C^{2}$ volume
preserving weakly hyperbolic group actions on closed manifolds.
For the integral action generated by a single Anosov
diffeomorphism this theorem is classical and originally due to
Anosov.

Motivated by the Franks/Manning classification of Anosov
diffeomorphisms on tori,
  we restrict our attention to weakly hyperbolic actions on the torus.  When the
acting group is a lattice subgroup of a semisimple Lie group with
no compact factors and all (almost) simple factors of real rank at
least two, we show that weak hyperbolicity in the original action
implies weak hyperbolicity for the induced action on the
fundamental group.  As a corollary, we obtain that any such action
on the torus is continuously semiconjugate to the affine action
coming from the fundamental group via a map unique in the homotopy
class of the identity. Under the additional assumption that some
partially hyperbolic group element has quasi-isometrically
embedded lifts of unstable leaves to the universal cover, we
obtain a conjugacy, resulting in a continuous classification for
these actions.
\end{abstract}

\section{Introduction}
In this article we investigate the notion of weak hyperbolicity
for group actions first introduced in \cite{MQ01}.   For linear
representations, weak hyperbolicity requires that there are no
nontrivial subrepresentations for which all eigenvalues of all
group elements have modulus one.  Weak hyperbolicity for a group
action on a closed manifold is a differential-geometric version of
that for representations.  The classical weakly hyperbolic
dynamical systems, i.e.  weakly hyperbolic actions of the
integers, correspond to the well understood class of Anosov
diffeomorphisms.  With this in mind, we shall draw conclusions
about weakly hyperbolic actions analogous to well known results
describing Anosov diffeomorphisms. For example, Theorem 3.5
establishes the ergodicity of weakly hyperbolic actions,
generalizing Anosov's work on ergodicity of Anosov diffeomorphisms
\cite[Theorem 4]{A}.  The proof of Theorem 3.5  uses ideas from
the well known Hopf argument (see \cite[Section 2.1]{BPSW} for a
description). Indeed, the proof is geometrically based on the
presence (and accessibility) of stable foliations, and technically
depends on the absolute continuity of these foliations.  Perhaps
new here is the fact that the core of this argument for ergodicity
lies in the use of a regularity result from \cite{RT} relating
Sobolev classes of functions measured tangentially with respect to
absolutely continuous foliations to global Sobolev classes of
functions.

The remainder of this work is motivated by that of Franks and
Manning (\cite{F70} , \cite{M74}) on the classification of Anosov
diffeomorphisms on tori. One may ask more generally whether weakly
hyperbolic actions on tori are classifiable. When the acting group
is a higher rank lattice, this question falls into Zimmer's
program for classifying volume preserving ergodic actions of
higher rank lattices on closed manifolds \cite{Z86}.  Our Theorem
5.2 is an analogue of Manning's contribution to the classification
of Anosov diffeomorphisms on tori. Roughly speaking, it asserts
that weak hyperbolicity  is inherited by the action on the
fundamental group when the acting group is a higher rank lattice.
The proof uses a rigidity property of the higher rank lattice,
specifically that such groups have Kazhdan's property (T), to
first draw conclusions about the action in the measurable category
and then uses the dynamical assumption of weak hyperbolicity to
bootstrap regularity from measurable to continuous.  Theorem 5.2
confirms what is suggested to be true by Margulis and Qian in
\cite{MQ01}.  Therein, and more generally in \cite{FW01}, the
analogue of Frank's contribution is proven in the higher rank
lattice setting.  Their result coupled together with the
complementary Theorem 5.2 establishes that all $C^2$ volume
preserving weakly hyperbolic actions of higher rank lattices on
tori, covered by an action on $\mathbb{R}^n$, are continuously
semiconjugate to the affine action coming from the fundamental
group (after possibly passing to a finite index subgroup of the
lattice). If, in addition, some group element acts by a partially
hyperbolic diffeomorphism with a quasi-isometric (in the universal
cover) unstable foliation, then we show this semiconjugacy is
injective, providing a continuous classification for this class of
actions. We believe that the continuous semiconjugacy is injective
without additional hypotheses and also believe a smooth
classification should be possible.

We are glad to thank Jeffrey Rauch for providing Theorem 2.3 and
outlining its Corollary 2.4, our advisor Ralf Spatzier, for his
ongoing help and encouragement, and the referee for invaluable
criticism on an earlier draft of this paper.

\section{Background}
This section contains the definitions and results with relevance
in subsequent sections.

\subsection{Partially Hyperbolic Dynamics}
Let $M$ be a closed Riemannian manifold and  let $\mu$ denote the
probability measure obtained by normalizing the Riemannian volume on $M$.
An element $f\in
\Diff^{k}(M)$ ($k\geq 1$) is a \textit{partially hyperbolic
diffeomorphism} if there exist continuous $df$-invariant
subbundles $E_{f}^{s},E_{f}^{c},E_{f}^{u}\subset TM$, and real numbers $C \ge
1, a > b \ge 1$ such that:

\begin{itemize}
\item $TM=E_{f}^{s}\oplus E_{f}^{c}\oplus E_{f}^{u}$ and

\item $\| d(f^{n})v^{u}\| \ge C^{-1}a^{n}\| v^{u}\|$, $\|
d(f^{n})v^{s}\| \le Ca^{-n}\|
v^{s}\|$, and  \\
$C^{-1}b^{-n}\| v^{c} \| \le \| d(f^{n})v^{c}\| \le Cb^{n}\|
v^{c}\|$, for all $v^{u}\in E_{f}^{u}$, $v^{c}\in E_{f}^{c}$,
$v^{s}\in E_{f}^{s}$ and all positive integers n.
\end{itemize}

We denote the set of $C^{k}$ volume preserving partially
hyperbolic diffeomorphisms by $PH_{\mu}^{k}(M)$.  The
distributions $E_{f}^{s}$ and $E_{f}^{u}$ are called the strong
stable and strong unstable distributions. It is well known (see
\cite[Theorem 5.5]{HPS} or \cite{P}) that these distributions
integrate uniquely and that the resulting integral submanifolds
are $C^{k}$ and form continuous foliations of $M$. These
foliations will be denoted by $\mathcal{W}_{f}^{s}$ and
$\mathcal{W}_{f}^{u}$ and are called the stable and unstable
foliations.  Recall that a
 \textit{d-dimensional continuous foliation} $\mathcal{W}$ of $M^{n}$ by
$C^{k}$ \textit{leaves} is a partition of $M$ into locally
immersed $C^{k}$ $d$-dimensional submanifolds called leaves so
that each point $x\in M$ has a \textit{foliated neighborhood}, i.e
a map $\Gamma:B^{d} \times B^{n-d} \rightarrow M$ where $B^{d}$
denotes the ball of dimension $d$ such that
\begin{itemize}
\item $\Gamma$ is a homeomorphism onto an open set in $M$ taking
(0,0) to $x$, \item for each $y\in B^{n-d}$, the map
$\Gamma(\cdot,y):B^{d}\rightarrow M$ belongs to $C^{k}(B^{d},M)$
and locally defines a leaf of the foliation, and \item the map
$B^{n-d}\rightarrow C^{k}(B^{d},M)$ given by $y \mapsto
\Gamma(\cdot,y)$ is continuous in the $C^{k}$ topology.
\end{itemize}  The leaf through the point $x$ is denoted by
$\mathcal{W}$(x).

To mimic the classical Fubini theorem, one would like that the
volume of a measurable set $E\subset M$ is obtained by integrating
the volume of $E$ in leaves along a transversal to the foliation.
A foliation $\mathcal{W}$ of $M$ is said to be \textit{absolutely
continuous} if for each open set $U$ of $M$ which is a union of
local leaves and each local transversal $T$ to the foliation there
is a measurable family of positive measurable functions
$\delta_x:\mathcal{W}(x)\cap U \rightarrow \mathbb{R}$ so that for
each measurable subset $E\subset U$,
$$\mu(E)=\int_{T} \int_{W(x)\cap U}
\chi_{E}(x,y)\delta_{x}(y)\,d\mu_{\mathcal{W}(x)}(y)\,d\mu_{T}(x),$$
where $\mu_{\mathcal{W}(x)}$ and $\mu_{T}$ denote the Riemannian
volumes in $\mathcal{W}(x)$ and $T$ respectively.  Absolute
continuity of a foliation as formulated above implies that zero
volume subsets in the foliated manifold have zero volume in leaves
through almost all points (\cite[Lemma 5.4]{B1}).  Surprisingly,
there are examples of foliations that are not absolutely
continuous (\cite{Mi}). A strictly stronger notion
(\cite[Proposition 3.5]{B1}) is that of \textit{transversal
absolute continuity} of a foliation. To define this, first note
that for any points $x_1 \in M$ and $x_2 \in \mathcal{W}(x_1)$ and
choice of transversals $T_i$ to the foliation through the points
$x_i$ ($i=1,2$), there is an associated \textit{Poincare map}.
This map is a homeomorphism $p:U_1 \rightarrow U_2$ between
neighborhoods $U_i$ of $x_i$ in $T_i$ satisfying $p(x_1)=x_2$ and
$p(x) \in \mathcal{W}(x)$ for each $x \in U_1$. A foliation
$\mathcal{W}$ is \textit{transversally absolutely continuous} if
all its Poincare maps are absolutely continuous maps with respect
to the induced Riemannian measures on the transversals. In other
words, for each choice of transversals $L_1$ and $L_2$ and
associated Poincare map $p$, there is a positive measurable
Jacobian $J:U_1 \rightarrow \mathbb{R}$ such that for each
measurable subset $A \subset U_1$,
$$\mu_{T_2}(p(A))=\int_{U_1} \chi_A(x)J(x) d\mu_{T_1}(x).$$
If, in addition these Jacobians are continuous and positive, then
the foliation is said to be \textit{measurewise} $C^1$. Of
technical importance is the following:

\begin{Thm}\cite[Theorem 2.1]{PS}
Let $f\in PH_{\mu}^{k}(M)$ and suppose $k\ge 2$.  Then the stable
and unstable foliations of $f$ are measurewise $C^1$.
\end{Thm}

In the course of the proof of Theorem 2.1, Pugh and Shub gave an
asymptotic expression for the Jacobians of the Poincare maps.
Starting from this expression, Nitica and Torok further
strengthened the absolute continuity property for stable and
unstable foliations.  They studied the regularity properties of
these foliations, first translating Pugh and Shub's work into a
statement about the absolute continuity of the local coordinate
charts for the foliation.  Their formulation shows that Jacobians
are differentiable along leaves and is useful for studying the
regularity properties of functions on $M$ that restrict to regular
functions on the leaves.  The following follows from the proof of
their regularity theorem:

\begin{Thm}\cite[Theorem 6.4]{NT}
Let $f\in PH^{2}_{\mu}(M)$.  Then there is a strong stable
foliation chart $\Gamma:B^{d}\times B^{n-d} \rightarrow M$ around
each $p\in M$ so that $\Gamma^{*}(\mu)=J(x,y)\,dx\,dy$ where the
Jacobian $J$ is an everywhere positive and continuous function on
$B^{d}\times B^{n-d}$ that has continuous (in $x$ and $y$
variables) partials of first order in the $x$ variables.
\end{Thm}

When a foliation $\mathcal{F}$ has local foliation charts as
described in Theorem 2.2, we shall say the foliation $\mathcal{F}$
is \textit{strongly absolutely continuous}.  The next theorem is a
local regularity result. It is used in subsequent sections to
reduce proving continuity of a function on a manifold to proving
that the function is Lipschitz when restricted to the stable
leaves of a family of partially hyperbolic diffeomorphisms whose
stable distributions jointly span the tangent bundle. Related
results include \cite{HuK}, \cite{J}, \cite{L},
 and \cite{LMM}.  First we establish some notation. For a
function $u$ on $M$, a strongly absolutely continuous foliation
$\mathcal{F}$, and $p\in (1,\infty)$, we say $u\in
H_{\mathcal{F}}^{1,p}(M)$ provided that for each foliation chart
as described by Theorem 2.2, $$\partial_{x_i}(u\circ \Gamma) \in
L^{p}(B^{d} \times B^{n-d}), i=1,\ldots d,$$ where differentiation
 is taken in the distributional sense.  Recently, Rauch and Taylor
 have established that functions in $H_{\mathcal{F}}^{1,p}(M)$
 are microlocally in the standard Sobolev space $H^{1,p}(M)$ away from the conormal
 bundle of $\mathcal{F}$ (\cite[Theorem 1.2]{RT}).  This implies the
 following:

\begin{Thm}\cite[Theorem 1.1]{RT}
Let $\mathcal{F}_{1},\ldots \mathcal{F}_{N}$ be strongly
absolutely continuous foliations of $M$.  Assume that for each $x
\in M$,
$$T_{x}M=\sum_{j=1}^{N} T_{x}\mathcal{F}_{j}.$$  Then, given $p\in
(1,\infty)$, if $u\in H_{\mathcal{F}_j}^{1,p}(M)$ for each $j\in
\{1,\ldots N\}$, then $u\in H^{1,p}(M)$, the standard Sobolev
space.
\end{Thm}

\begin{Cor}
Let $\mathcal{F}_{1},\ldots \mathcal{F}_{N}$ be strongly
absolutely continuous foliations of $M$. Assume that for each
$x\in M$,
$$T_{x}M=\sum_{j=1}^{N} T_{x}\mathcal{F}_{j}.$$ Suppose that $u\in
L^{2}(M)$ and that there is a constant $K>0$ so that the
restrictions of $u$ to almost all leaves of the foliations
$\{\mathcal{F}_{j}\}$ are almost everywhere $K$-Lipshitz, then $u$
agrees almost everywhere with a continuous function.
\end{Cor}

\begin{Prf}
As a first step we argue that $u$ has bounded tangential first
order derivatives parallel to the leaves of the foliations. These
derivatives are only taken in the distributional sense. In a
second step, we will use regularity bootstraping techniques to
finish the proof. To this end, fix $j \in \{1, \ldots N \}$, $p
\in M$, and let $d=\dim \mathcal{F}_j$. As the regularity
properties of $u$ are a local matter, we may identify a local
paramaterization of $\mathcal{F}_j$ in a neighborhood of $p \in M$
with
 a foliated open neighborhood $\Omega \subset \mathbb{R}^{n}$, described by
  a paramaterization $$\Gamma: U \times V \rightarrow \Omega,$$ with
   $U \subset \mathbb{R}^d$,$V \subset \mathbb{R}^{n-d}$ and Jacobian
    $J(x,y)$ as in Theorem 2.2.   Let $X$ be a vector field in $\Omega$
     such that $X\circ \Gamma = \sum_{j=1}^{d} a_{j}(x,y) \partial_{x_j}$,
      with $\partial_{x_i}a_{j}(x,y) \in C^{0}(U \times V)$ for
       $i=1,\ldots,d$ (such a field is well defined, independent of
        the choice of ``flattening" paramaterization).  Let
         $U' \subset \subset U$, $V' \subset \subset V$, denote by
          $\mathcal{D}^\prime(U' \times V')$
the space of distributions on $U' \times V'$, and define
 $$X: L^{1}_{loc}(U' \times V') \rightarrow \mathcal{D}^\prime(U' \times V'),$$
 by $X(f)(\phi)=\int_{V'}\int_{U'} f (-\sum_{j=1}^{d} \partial_{x_j}(\phi a_j))dx \, dy$
  for $f\in L^1_{loc}(U' \times V')$ and $\phi \in C_{c}^{\infty}(U' \times V')$.
   Since $M$ has finite volume, $u\in L^{1}(M)$, and by continuity and positivity of
   the Jacobian $J$, it follows that $u\circ \Gamma \in L^{1}_{loc}(U' \times V')$.
   To complete this step, we must show that
   $X(u\circ \Gamma) \in L^{\infty}(U' \times V') \subset \mathcal{D}^\prime(U' \times V')$.
    To this end we define (when n is sufficiently large) the difference
     quotients $u_{n,j}(x,y):= n(u\circ \Gamma(x+\frac{1}{n}e_j,y) - u\circ \Gamma(x,y))$
     and the functions $g_{n}:= \sum_{j=1}^{d} u_{n,j}a_j$.  By the hypothesis on $u$ and
     the continuity of the $a_j$, each $g_{n} \in L^{\infty}(U' \times V')$ and there is a
      uniform $C > 0$ so that $||g_{n}||_{\infty} < C$.  In particular there is a
       uniform bound $C'$ on the $L^{2}$ norms of the $g_n$, so that a subsequence of the
        $g_n$ weakly converge to a function $g \in L^{2}(U' \times V')$.
This weak limit is also a distributional limit and furthermore it follows from
 the uniform bound $||g_n||_{\infty} < C$ that $||g||_{\infty} < C$.   On the
  other hand, a routine change of variables argument shows that the $g_n$ converges
   to $X(u \circ \Gamma)$ in $\mathcal{D}^\prime
(U' \times V')$, completing this step.

By the first step and since $u\in L^{2}(M)$,  $u\in H^{1,2}_{\mathcal{F}_j}(M^{n})$ for
  each $j \in \{1,\ldots, N\}$.  Applying Theorem 2.3 establishes that
   $u\in H^{1,2}(M^{n})$.  By the Sobolev imbedding theorem (see \cite[Proposition 2.2]{T}),
   $H^{1,p}(M^{n}) \subset L^{\frac{np}{n-p}}(M^n)$
   for $p\in [1,n)$. From this, it follows that $u \in L^{\frac{2n}{n-2}}(M^{n})$.
    Again by the first step, we now obtain $u\in  H^{1,\frac{2n}{n-2}}_{\mathcal{F}_j}(M^{n})$
    for each $j$ and hence $u \in H^{1, \frac{2n}{n-2}}(M^{n})$ by Theorem
    2.3.  Note that for all $0<y<n$, $\frac{ny}{n-y}-y$ is positive and increasing
     as a function of $y$.
     We may therefore define a finite increasing sequence $\{x_0, \ldots x_J\}$
      inductively by letting $x_0=2$ and $x_i=\frac{nx_{i-1}}{n-x_{i-1}}$
       whenever $x_{i-1}<n$. Repeating the above argument shows that $u\in
       H^{1,x_J}(M^n)$. When $x_J=n$, note that $u\in H^{1,n-\epsilon}(M)$ for small
         $\epsilon>0$.  Choose $\epsilon< \frac{n}{2}$ so that
$\frac{n(n-\epsilon)}{n-(n-\epsilon)}>n$.  Then
         $u\in H^{1,\frac{n(n-\epsilon)}{n-(n-\epsilon)}}(M)$.  Therefore,
         $u\in H^{1,p}(M^n)$ for some $p>n$ whether $x_J>n$ or $x_J=n$. The conclusion
           of the corollary follows by the Sobolev
           imbedding (see \cite[Proposition 2.4]{T}
 $H^{1,p}(M^{n}) \subset C^{0}(M^{n})$ for
 $p>n$.
\hfill $\Box$
\end{Prf}

\subsection{Kazhdan's property T}

Let $\Gamma$ be a discrete group and
$\mathcal{H}$ a separable Hilbert space. A unitary
representation $\pi:\Gamma \rightarrow U(\mathcal{H})$  is said to have
\textit{almost invariant vectors} if and only if for every
$\epsilon
> 0$ and each compact subset $K\subset \Gamma$, there exists a unit
vector $u(\epsilon,K) \in H$ so that $\| \pi(\gamma)u-u \| <
\epsilon$ for each $\gamma\in K$.  The group $\Gamma$ has
\textit{Kazhdan's property (T)} if every unitary representation of
$\Gamma$ with almost invariant vectors has a nonzero invariant
vector. It is well known that higher rank lattices have this fixed
point property.  There is a useful characterization of this
property in terms of the vanishing of the first cohomology group
of $\Gamma$ with coeficients in a unitary representation.  Let
$\phi$ be a representation of $\Gamma$ on a vector space $V$.  A
\textit{1-cocycle} is a map $c:\Gamma \rightarrow V$ such that
$c(gh)=c(g) + \phi(g)(c(h))$ for each $g,h \in \Gamma$. A
\textit{coboundary} is a cocycle $c_{v}$ of the form
$c_{v}(g)=v-\phi(g)(v)$ for some fixed $v\in V$.  The first
cohomology of $\Gamma$ with coeficients in the representation
$\phi$ is the group of cocycles modulo the normal subgroup of
coboundaries and is denoted by $H^{1}(\Gamma,\phi)$.

\begin{Thm}\cite{D}, \cite{G}
A group $\Gamma$ has Kazhdan's property if and only if
$H^{1}(\Gamma,\phi)=0$ for every unitary representation $\phi$ of $\Gamma$.
\end{Thm}

\section{Weak Hyperbolicity, Accessibility, and Ergodicity}
The purpose of this section is to define the actions under
consideration as well as the associated notion of accessibility.
Recall that stable-unstable accessibility plays a central role in
Hopf's argument for proving the ergodicity of Anosov systems. Here
too, accessibility enters into a proof that $C^{2}$ weakly
hyperbolic actions of discrete groups on closed manifolds are
ergodic. In what follows, let $M$ denote a closed Riemannian
manifold.

\begin{Def}
A \textit{weakly hyperbolic family on M} is a finite subset of
$PH_{\mu}^{1}(M)$ which is infinitesimally accessible in stable
directions. More precisely, it is a finite family
$\{\gamma_{1},\ldots, \gamma_{k}\} \subset PH_{\mu}^{1}(M)$ with
associated splittings $TM=E_{i}^{s}\oplus E_{i}^{c} \oplus
E_{i}^{u}$ such that $TM=\sum_{i=1}^{k} E_{i}^{s}$.  An action of
a discrete group $\Gamma$ on $M$ is \textit{weakly hyperbolic} if
there is a finite family of group elements that act as a weakly
hyperbolic family on $M$.
\end{Def}
The existence of the stable
foliations suggests that infinitesimal accessibility along stable
directions integrates to a notion of local accessibility along
stable leaves.

\begin{Def} Let $F=\{\gamma_1,\ldots, \gamma_k\}$ be a finite
family of partially hyperbolic diffeomorphims on $M$ with
associated stable foliations $\{\mathcal{W}^{s}_{1},\ldots,
\mathcal{W}^{s}_{k}\}$. An \textit{admissible path} for $F$ is a
path in $M$ that is piecewise sequentially contained in leaves of
the stable foliations: i.e. a path $a:[0,T] \rightarrow M$ with a
subdivision $t_0=0 \leq t_1 \leq \cdots \leq t_k=T$ such that
$a([t_{j-1},t_{j}]) \subset \mathcal{W}_{j}^{s}(a(t_{j-1}))$ for
each $j \in \{1, \cdots, k\}$. The family  $F$ is \textit{locally
accessible} if for each $x\in M$ there is an open neighborhood
$U_{x}$ of $x$ such that each $y \in U_{x}$ is the endpoint of an
admissible path for $F$ beginning at $x$.
\end{Def}

\begin{Lem}
A weakly hyperbolic family $F$ of partially hyperbolic
diffeomorphisms on a closed manifold $M$ has the local
accessibility property.
\end{Lem}

\begin{Prf}
Let $F=\{\gamma_1, \ldots, \gamma_k\}$ be a weakly hyperbolic
family on $M$ and fix $x\in M$.  Note that any piecewise
differentiable path in $M$ beginning at $x$ which is piecewise
sequentially tangent to the stable distributions
$\{E_{i}\}_{i=1}^{k}$ is an admissible path for $F$ since the
stable distributions uniquely integrate.  As this is a local
question, we may view each stable distribution $E_{i}$ near $x$ as
being given by the span of a finite ordered family of continuous
vector fields $\{X_{i}^{1},\ldots,X_{i}^{d(i)}\}$, where
$d(i)=\dim(E_{i})$.  Any curve in $M$ which sequentially defines a
solution to the fields $\{X_{i}^{j}\}_{j=1}^{d(i)}$ is tangent to
$E_{i}$. Therefore, any curve that sequentially defines a solution
to the combined ordered family of fields $X=\{X_{1}^{1}, \ldots,
X_{1}^{d(1)}, \ldots , X_{k}^{1}, \ldots, X_{k}^{d(k)}\}$ is
admissible for the family $F$.  The hypothesis that the family $F$
is weakly hyperbolic implies that the family $X$ locally spans the
tangent bundle near the point $x$.  Establishing local
accessibility for a weakly hyperbolic family therefore reduces to
the next lemma.
\end{Prf}

\begin{Lem}
Let $\{X_1,\ldots, X_N\}$ be a family of nonvanishing continuous
vector fields in $\mathbb{R}^n$ such that $\spn(X_1,\ldots,
X_N)=\mathbb{R}^n$. Then the set of endpoints of curves beginning
at $0$ which sequentially define solutions to these fields
contains an open set around $0$. \hfill $\Box$
\end{Lem}

In \cite{KP} they study control systems with low regularity.  Their
Corollary 4.5 implies this last lemma.  Note that if the vector fields
had enough regularity to have $C^1$ exponential maps,
then the spanning hypothesis implies that the composition of the
exponential maps is a submersion at $0 \in \mathbb{R}^n$, from
which the lemma follows.  For continuous fields, solutions need
not be unique and consequently, there are no exponential maps.  The
result in \cite{KP} works in a significantly more general context and
circumvents non-uniqueness by perturbing a system with low
regularity to one with better regularity and applying a degree
theory for set valued maps.  An alternative approach is to
approximate the family of continuous fields by a sequence of
families of smooth fields and then argue that the sequence of
open sets one obtains for the approximating families have a
uniform lower bound on their size in order to take limits.  For
brevity, we omit the details of this direct approach.

In the remainder of this section we use local accessibility
together with the regularity result, Corollary 2.4, in order to
prove that $C^{2}$ weakly hyperbolic volume preserving actions of
discrete groups are ergodic.  Let $\Gamma$ be a discrete group.
Recall that proving ergodicity for a measure preserving $\Gamma$
action on a Borel probability space $(X,\mu)$ is equivalent to
showing that any square integrable function which is almost
everywhere invariant under the action of $\Gamma$ is almost
everywhere constant. To accomplish this, we first show that a
square integrable function almost everywhere invariant under an
element $\gamma \in PH_{\mu}^{2}(M)$ has the needed tangential
regularity in order to apply Corollary 2.4.  In fact, such a function
will be essentially constant on almost all leaves
of $\mathcal{W}_{\gamma}^{s}$. Weak hyperbolicity together with
Corollary 2.4 then implies that each almost
everywhere invariant representative of an $L^{2}$ element is
necessarily equivalent to a continuous function $f$ which is
almost everywhere invariant.  As $f$ is continuous, it is
everywhere invariant and hence constant on stable leaves. Local
accessibility implies that $f$ is locally constant and therefore
constant on $M$.

What follows makes this reasoning precise.  Recall that for each
volume preserving transformation $T$ of $M$, $\mu$ has an ergodic
decomposition (\cite[Theorem 2.19]{FMW}); more precisely, there is
a full volume $T$-invariant Borel set $M'\subset M$, a standard
Borel probability space $(\Omega,\nu)$, a Borel map $\xi:\Omega
\rightarrow \Prob(M')$, and a $T$-invariant Borel map $\psi:M'
\rightarrow \Omega$ such that $\xi(\omega)(\psi^{-1}(\omega))=1$
for each $\omega \in \Omega$, $\mu=\int_{\Omega}
\xi(\omega)\,d\nu(\omega)$, and $\xi(\omega)$ is quasi-invariant
and ergodic for each $\omega \in \Omega$, where $\Prob(M')$
inherits its Borel structure from the weak* topology.

\begin{Lem}
Suppose that $T \in PH^{1}_{\mu}(M)$ has an absolutely continuous
stable foliation and that a square integrable function $f$ is
almost everywhere $T$-invariant. Then for almost all leaves, $f$
restricted to the leaf is almost everywhere constant.
\end{Lem}

\begin{Prf}
For $g\in L^{1}(M,\mu)$ we denote the set of Birkhoff regular
points for $g$ by $B_{g}:=\{x \in M' | g^{+}(x)=\int_{M'} g\,
d\xi(\psi(x))\}$, where $g^{+}(x)$ is given by
 $\lim_{n \rightarrow \infty}\frac{1}{n} \sum_{i=0}^{n-1} g(T^{i}(x))$
 at points where this
 limit converges.  It follows easily from the description of the
ergodic decomposition that these sets are Borel sets and the
Birkhoff ergodic theorem implies these sets have full volume in
$M$. From the separability of $C^{0}(M)$, it follows that the set
$B_{0}:=\cap_{g\in C^{0}(M)} B_{g}$ has full volume in $M$.  By
hypothesis there is a full volume $T$-invariant set $I\subset M$
so that $f$ is $T$-invariant in $I$.  By absolute continuity of
the stable foliation, there is a full volume subset of good points
$G\subset (I\cap B_{f} \cap B_{0})$ such that $x\in G$ implies
$(I\cap B_{f} \cap B_{0})$ is conull in $\mathcal{W}^{s}_{T}(x)$.
Fix $x\in G$ and let $y_{1}, y_{2} \in \mathcal{W}^{s}_{T}(x) \cap
(I \cap B_{f} \cap B_{0})$.  We first argue that $y_1$ and $y_2$
lie in the same ergodic component, or more precisely, that
$\xi(\psi(y_{1}))=\xi(\psi(y_{2}))$.  Indeed the Hopf argument
(see e.g. \cite[Section 2.1]{BPSW}) shows that for any continuous
function $g$, whenever $g^{+}(x):= \lim_{n \rightarrow
\infty}\frac{1}{n} \sum_{i=0}^{n-1} g(T^{i}(x))$ converges,
$g^{+}(y)$ converges to $g^{+}(x)$ for all $y\in \mathcal
W_{T}^{s}(x)$. Since $y_{1}, y_{2} \in B_{0}$, it follows that
$\int_{M} g\, d \xi(\psi(y_1)) = \int_{M} g\, d \xi(\psi(y_2))$
for all continuous functions $g$, whence
$\xi(\psi(y_{1}))=\xi(\psi(y_{2}))$.  Since $y_{1},y_{2} \in
(I\cap B_{f})$,
$$f(y_{1})=f^{+}(y_{1})=\int_{M'}f \,d\xi(\psi(y_{1}))=\int_{M'}f\,
d\xi(\psi(y_{2}))=f^{+}(y_{2})=f(y_2).$$

\hfill$\Box$
\end{Prf}

From Corollary 2.4 and Lemma 3.3 we deduce the following:
\begin{Cor}
Suppose that $f$ is a square integrable function on $M$ that is
almost everywhere invariant under a weakly hyperbolic $C^{2}$
family on $M$.  Then $f$ is almost everywhere equal to a
continuous function $g$.
\end{Cor}

\begin{Thm}
Let $\rho$ be a $C^{2}$ volume preserving weakly hyperbolic action
of a discrete group $\Gamma$ on $M$. Then the $\Gamma$ action is ergodic.
\end{Thm}

\begin{Prf}
Let $f$ be any square integrable almost everywhere
$\Gamma$-invariant function.  By the last corollary, $f$ is almost
everywhere equal to a continuous function $g$ that is almost
everywhere $\Gamma$-invariant. By continuity of $g$ and since full
volume sets are dense, $g$ is everywhere invariant. Hence $g$ is
constant on all stable leaves of all the elements in the weakly
hyperbolic family. By local accessibility, $g$ is locally constant
and hence constant.  Therefore $f$ is essentially constant.
\hfill$\Box$
\end{Prf}

\section{Weak Hyperbolicity is Inherited}
In this section we consider volume preserving weakly hyperbolic
actions of discrete property (T) groups $\Gamma$ on tori. We argue
that whenever such an action is covered by a $\Gamma$ action on
$\mathbb{R}^n$, the representation coming from the homomorphism
$\Gamma \rightarrow \Out(\pi_1(\mathbb{T}^n))$ cannot split as a
nontrivial direct sum of subrepresentations, one of which is
isometric. To this end, let $\Gamma$ denote a discrete Kazhdan
group and $\rho:\Gamma \rightarrow \Diff^{2}(\mathbb{T}^n)$ be a
volume preserving weakly hyperbolic action covered by an action,
$\overline{\rho}:\Gamma \rightarrow \Diff^2(\mathbb{R}^n)$.  The
action $\rho$ induces a homomorphism $\Gamma \rightarrow
\Out(\pi_1(\mathbb{T}^n))$ that lifts to a homomorphism
$\pi:\Gamma \rightarrow \Aut(\mathbb{Z}^n)$. Note that $\pi$
satisfies
$$\overline{\rho}(\gamma)(x+z)=\overline{\rho}(\gamma)(x) + \pi(\gamma)(z),$$
for each $\gamma \in \Gamma$, $x \in \mathbb{R}^n$, and $z \in
\mathbb{Z}^n$.  For simplicity, we also let $\pi$ denote the
representation $\Gamma \rightarrow \GL(n,\mathbb{R}^n)$ induced by
the homomorphism $\pi$.  Finally, let
 $\rho_{0}$ denote the linear action
on the torus induced by $\pi$.  We assume that there is a direct
sum decomposition of $\mathbb{R}^n$ into $\Gamma$ invariant
subspaces $C$ and $H$ so that the restriction of $\pi$ to $C$,
$\pi^C$, acts isometrically in $C$ and will argue that $C=\{0\}$.
The main idea is that the representation $\pi$ coarsely
approximates
 the lifted action $\overline{\rho}$.  Since all tangent directions are spanned by directions
  that are contracted uniformly by some element of the group under $\overline{\rho},$ the
  approximating action $\pi$ cannot have a nontrivial isometric invariant subspace $C$.
  Making this line of reasoning precise involves analyzing a
   cocycle (first introduced in \cite{MQ01}) which measures the difference between the
   lifted action and the induced action on the fundamental group.  The heart of the argument
    lies in the following:

\begin{Prop} Assume the hypotheses above.  Then there is a
continuous map $\phi:\mathbb{R}^n \rightarrow C$ of the form
$\phi(x)= \proj_{C}(x)+\sigma([x])$ where $\sigma \in
C^{0}(\mathbb{T}^n,C)$ such that
$\pi^{C}(\gamma)(\phi(x))=\phi(\overline{\rho}(\gamma)x)$.
\end{Prop}

\begin{Prf}
Define the map $A:\Gamma \times \mathbb{R}^n \rightarrow
\mathbb{R}^n$ by the equation
$$\overline{\rho}(\gamma)x=\pi(\gamma)(x+ A(\gamma,x)).$$  In view of how $\pi$ is defined,
 for a fixed $\gamma\in \Gamma$, the
function $A(\gamma,\cdot)$ descends to a function on the torus.
 Following $A$ by the projection to $C$ parallel to the
complementary subspace
 $H$, we obtain a
map $A^{C}:\Gamma \times \mathbb{T}^n \rightarrow C$.  Since $H$
is an invariant complement to $C$, $A^C$ solves the equation
$$\proj_{C}(\overline{\rho}(\gamma)x)=\pi^{C}(\gamma)(\proj_{C}(x)+A^{C}(\gamma,[x])).$$
Define the unitary representation $\lambda:\Gamma \rightarrow
U(L^{2}(\mathbb{T}^n,C))$ by
$$(\lambda(\gamma)f)([x]):=\pi^C(\gamma)f(\rho(\gamma^{-1})[x]).$$

Next, we check that the map $c:\Gamma \rightarrow
L^{2}(\mathbb{T}^n,C)$ given by $\gamma \mapsto
A^{C}(\gamma^{-1})$ is a  1-cocycle in $Z^{1}(\Gamma,
\lambda)$.  As a preliminary step, we establish
 $$(*)\>  A^{C}(\gamma_1 \gamma_2,[x])=A^{C}(\gamma_2,[x])+
 \pi^{C}(\gamma_2^{-1})A^{C}(\gamma_1,\rho(\gamma_2)[x]). $$  Indeed,  \begin{align*}
\pi(\gamma_1 \gamma_2)(x+ A(\gamma_1 \gamma_2,[x]))
&=& \overline{\rho}(\gamma_1)(\overline{\rho}(\gamma_2)x) &=&
 \pi(\gamma_1)(\overline{\rho}(\gamma_2)x + A(\gamma_1,\rho(\gamma_2)[x]))=
\end{align*}
$$\pi(\gamma_1)[ \pi(\gamma_2)(x+ A(\gamma_2,[x])) +
 A(\gamma_1,\rho(\gamma_2)[x])] =$$
$$ \pi(\gamma_1 \gamma_2)(x + A(\gamma_2,[x]) +
\pi(\gamma_2^{-1})A(\gamma_1,\rho(\gamma_2)[x])). $$
Applying $\pi(\gamma_2^{-1} \gamma_1^{-1})$ to the first and last term and projecting to the
subspace $C$ gives $(*)$.  Therefore,
$$A^{C}(\gamma_2^{-1} \gamma_1^{-1},\cdot)=A^{C}(\gamma_1^{-1},\cdot)
+ \pi^{C}(\gamma_1)A^{C}(\gamma_2^{-1},\rho(\gamma_1^{-1})\cdot)=
A^{C}(\gamma_1^{-1},\cdot)+\lambda(\gamma_1)A^{C}(\gamma_2^{-1},\cdot),$$
showing $c \in Z^{1}(\Gamma,\lambda)$.

By Theorem 2.5 there is some $\sigma \in L^{2}(\mathbb{T}^n,C)$ such that
$$(**)\>  A^{C}(\gamma^{-1})= \sigma - \lambda(\gamma)\sigma,$$ holds as an
equation in $L^{2}$ for every $\gamma \in \Gamma$.  Next, we argue
that $\sigma$ agrees almost everywhere with a continuous function,
and consequently, that $(**)$ holds as an equation in
$C^{0}(\mathbb{T}^{n},C)$.  In view of Corollary 2.4, it suffices
to show that for a partially hyperbolic diffeomorphism
$\rho(\gamma) \in PH_{\mu}^{2}(\mathbb{T}^{n})$, there is a
uniform constant $C >0$ such that the restriction of $\sigma$ to
almost all leaves of the stable foliation
$\mathcal{W}^{s}_{\rho(\gamma)}$ is almost everywhere
$C$-Lipschitz.  To show this, we first simplify notation writing
$\gamma$ instead of $\rho(\gamma)$, $\hat{\gamma}$ instead of
$\pi^{C}(\gamma^{-1})$, and will let $f$ denote the $C^{2}$
function $A^{C}(\gamma)$.  In this simplified notation, $(**)$
implies that there is a full volume $\gamma$-invariant set $I
\subset \mathbb{T}^n$ so that $$(***)\>   \sigma(x)=f(x) +
{}^{\hat{\gamma}} \sigma({}^{\gamma} x)$$ for $x\in I$.  By the
ergodic decomposition (\cite[Theorem 2.19]{FMW}), there is a full
volume $\gamma$-invariant set $M'\subset \mathbb{T}^{n}$, a
standard Borel probability space $(\Omega,\nu)$, a
$\gamma$-invariant Borel map $\psi:M' \rightarrow \Omega$, and a
Borel map $\xi:\Omega \rightarrow \Prob(M')$ such that
$\mu=\int_{\Omega} \xi(\omega) d \nu$, where each $\xi(\omega)$ is
a quasi-invariant ergodic probability measure.  As in Lemma 3.3,
for each $g\in L^{1}$ define the full volume set
 $B_{g}\subset M'$ by  $B_{g}:=\{x \in M' | g^{+}(x)=\int_{M'} g\,
d\xi(\psi(x))\}$ and the full volume set $B_{0}= \cap_{g\in C^{0}}B_g$.  By
 Lusin's theorem there is a sequence of compact
sets $K_{j} \subset K_{j+1}$ such that the restriction of $\sigma$
to $K_j$ is uniformly continuous and $\mu(K_j)>1-\frac{1}{2^{j}}$.
Let $K=\bigcup_{j} K_{j}$ and $B_{K}:=\{x\in M'|
\xi(\psi(x))(K)=1\}$. It is straightforward to argue that $B_{K}$
is a full volume Borel set.  Finally, let $G:= I \cap B_{K}
\cap_{j} B_{\chi(K_j)} \cap B_{0}$.  Then $G$ is a full volume
Borel set and by absolute continuity, there is a full volume
subset of  points $E\subset G$ so that $G$ is conull in
$\mathcal{W}^{s}_{\gamma}(x)$ whenever $x\in E$. Fix $x\in E$ and
$y_1,y_2$ in $G\cap \mathcal{W}^{s}_{\gamma}(x)$. Since $y_1,y_2
\in B_{0}$, the Hopf argument (see Lemma 3.3) shows that
$\xi(\psi(y_1))=\xi(\psi(y_2))$.  Let $m$ denote this ergodic
measure. Since $y_1 \in B_{K}$, $m(K)=1$ so that there is a large
enough $j$ for which $m(K_j)>\frac{1}{2}$.  Since $y_1,y_2 \in
B_{\chi(K_j)}$ there are infinitely many $n$ for which
$\gamma^{n}y_1$ and $\gamma^{n} y_2$ both lie in $K_j$.  Let
$\epsilon >0$.  As $\gamma$ contracts distances in
$\mathcal{W}^{s}_{\gamma}(x)$ by some constant $\lambda<1$ and by
uniform continuity of $\sigma$ in $K_j$, there is a large enough
$N$ so that $d_{C}(\sigma(\gamma^{N} y_1), \sigma(\gamma^{N} y_2))
< \epsilon$.  By iterating $(***)$, we therefore obtain that
$$d_C(\sigma(y_1),\sigma(y_2))
=\|\sum_{i=0}^{N-1} \widehat{\gamma^{i}}(f(\gamma^{i} y_{1})-f(\gamma^{i}
y_{2})) +(\sigma(\gamma^{N} y_{1})-\sigma(\gamma^{N} y_{2}))\| $$
$$\le \sum_{i=0}^{N-1} L d_{M}(\gamma^{i} y_{1},\gamma^{i} y_{2})
+ \epsilon \le \sum_{i=0}^{\infty} L \lambda^{i}
d_{M}(y_{1},y_{2}) +\epsilon,$$ where $L$ is a Lipshitz constant
for $f$ and $\lambda$ is the contraction constant for $\gamma$.

Finally, we show that the map $\phi:\mathbb{R}^{n} \rightarrow C$
defined by $x \mapsto \proj_{C}(x) + \sigma([x])$ is equivariant:
$$\phi(\overline{\rho}(\gamma)x)=
 \proj_{C}(\overline{\rho}(\gamma)x) + \sigma(\rho(\gamma)[x])=\proj_{C}(\pi(\gamma)x+
 \pi(\gamma)A(\gamma,[x])) + \sigma(\rho(\gamma)[x])=$$ $$\pi^{C}(\gamma)\proj_{C}(x)
 +\pi_{C}(\gamma)A^{C}(\gamma,[x]) +\sigma(\rho(\gamma)[x])=\pi^{C}(\gamma)\proj_{C}(x)
  + \pi^{C}(\gamma)\sigma([x])=$$ $$\pi^{C}(\gamma)\phi(x),$$ where the penultimate equality
   follows from  $(**)$.

\hfill $\Box$
\end{Prf}

\begin{Lem}
With the same hypotheses as in Prop 4.1, the map $\phi$ is
constant.
\end{Lem}

\begin{Prf}
Fix a lift $\widetilde{\mathcal{W}}^{-}(x)$ of a stable leaf of a
partially hyperbolic diffeomorphism $\rho(\gamma)$ and let $y\in
\widetilde{\mathcal{W}}^{-}(x)$. Suppose that
$d_{0}:=d_{C}(\phi(x),\phi(y))
> 0.$ By uniform continuity of $\phi$ there is a
$\delta > 0$ so that $d(x,y) < \delta$ implies
$d_{C}(\sigma(x),\sigma(y)) < \frac{d_{0}}{2}$.  For sufficiently
large $n$,
$d(\overline{\rho}(\gamma^{n})x,\overline{\rho}(\gamma^{n})y)<
\delta$. Therefore,
$$d_{0}=d_{C}(\phi(x),\phi(y))=
d_{C}(\pi^{C}(\gamma^{n})\phi(x),\pi^{C}(\gamma^{n})\phi(y)))=
d_{C}(\phi(\overline{\rho}(\gamma^{n})x),\phi(\overline{\rho}(\gamma^{n})y))
< \frac{d_{0}}{2}.$$  This yields a contradiction unless $\phi$ is
constant on lifts of stable leaves and therefore locally constant
by the local accessibility property for weakly hyperbolic
families.

\hfill $\Box$
\end{Prf}

\begin{Thm}
With the same hypotheses as in Prop 4.1, the subspace $C$ is
trivial.
\end{Thm}

\begin{Prf}
Suppose not.  Since $\phi$ is defined by adding a bounded map to
the projection map to the subspace $C$, $\phi$ is unbounded,
contradicting Lemma 4.2. \hfill $\Box$
\end{Prf}

\section{Towards Global Rigidity}
In this section we assume that the acting group $\Gamma$ is a
lattice in a connected semisimple real Lie group $G$ with each
(almost) simple factor having real rank at least two.  For
simplicity we call such a group a higher rank lattice.  In this
section, we deduce that all weakly hyperbolic $C^{2}$ volume
preserving actions on a torus that lift to the universal cover are
semiconjugate to the linear action coming from the fundamental
group when restricted to a finite index subgroup of the acting
lattice. We also argue that this semiconjugacy is injective under
the additional hypothesis that the leaves of the lift of the
unstable foliation of a partially hyperbolic group element are
quasi-isometrically embedded in $\mathbb{R}^n$.

Let $\rho: \Gamma \rightarrow \Diff^{2}(\mathbb{T}^{n})$ be a
volume preserving action that lifts to an action on the universal
cover, $\overline{\rho}:\Gamma \rightarrow
\Diff^{2}(\mathbb{R}^n)$. Let $\pi_{\rho}: \Gamma \rightarrow
\GL(n,\mathbb{Z})$ be the induced homomorphism and $\rho_{0}$ be
the affine action induced by $\pi_{\rho}$. In \cite{MQ01}, they
make the following:

\begin{Def}
The representation $\pi_{\rho}:\Gamma \rightarrow
\GL(n,\mathbb{R})$ is said to be weakly hyperbolic if the Zariski
closure of $\pi(\Gamma)\subset \GL(n,\mathbb{R}^n)$ is not
precompact in any of its nontrivial subrepresentations.
\end{Def}

The next theorem is from \cite{MQ01}.  We remark that although
they assume the existence of a periodic point for the action,
their proof works equally well under the weaker assumption that
the action lifts to $\mathbb{R}^n$.
\begin{Thm}\cite[Theorem 6.10]{MQ01}
Let $\Gamma$, $\rho$, and $\overline{\rho}$ be as above.  If
$\pi_{\rho}$ is weakly hyperbolic, then there exists a finite
index subgroup $\Gamma^{'} < \Gamma$ and a map $\phi \in
C^{0}(\mathbb{T}^n)$, unique in the homotopy class of the identity
such that $\phi \circ \rho(\gamma)= \rho_{0}(\gamma)\circ \phi$
for all $\gamma \in \Gamma^{'}$.
\end{Thm}

Next we establish a complimentary statement.
\begin{Thm}
Let $\Gamma$, $\rho$, and $\overline{\rho}$ be as above.  If
$\rho$ is weakly hyperbolic, then $\pi_{\rho}$ is weakly
hyperbolic.
\end{Thm}

\begin{Prf}
Suppose that the Zariski closure of $\pi(\Gamma)$ is precompact in
$\GL(C)$ for some invariant subspace $C$.  Since the Zariski
closure of $\pi(\Gamma)$ is semisimple (\cite[Prop. IX
5.7]{Mar91}, there is an invariant subspace complementary to $C$.
Since $\Gamma$ has property (T), Theorem 4.1 implies $C={0}$.
\hfill $\Box$
\end{Prf}

\begin{Cor}
Let $\Gamma$, $\rho$, and $\overline{\rho}$ be as above.  If
$\rho$ is weakly hyperbolic, then after passing to a finite index
subgroup of $\Gamma$, $\rho$ is $C^{0}$-semiconjugate to the
affine action coming from the homomorphism $\Gamma \rightarrow
\Out(\pi_1(\mathbb{T}^n))$ by a map which is unique in the
homotopy class of the identity.
\end{Cor}

{\bf Remark}: For many lattices, the above results hold without
assuming that the lattice action $\rho$ on the torus is covered by
an action $\overline{\rho}$ of $\mathbb{R}^n$. Indeed, the action
will lift whenever $H^2(\Gamma, \mathbb{Z}^n)=0$ (see e.g.
\cite{FW01}). For relevant cohomology vanishing results see
\cite{Bo}, \cite{BoW}, \cite{Zu}, and \cite{K}.

\vskip 10pt

Next we argue that this semiconjugacy is a $C^{0}$ conjugacy
provided that the lifts of unstable leaves to $\mathbb{R}^n$ have
intrinsic distances comparable to Euclidean distance for some
partially hyperbolic group element.  Following \cite{B} we make the
following:

\begin{Def}
A foliation $\mathcal{W}$ of a simply connected metric space
$(X,d)$ is said to be \textit{quasi-isometric} if there are
uniform constants $a,b>0$ such that for each $x\in X$ and $y\in
\mathcal{W}(x)$, $d_{\mathcal{W}(x)}(x,y) \leq a d_{X}(x,y) + b$.
\end{Def}

\begin{Prop}
Let $\Gamma$, $\rho$, $\overline{\rho}$, $\pi$, and $\rho_{0}$ be
as above. Let $\phi$ be the unique continuous map homotopic to the
identity such that $\phi \circ \rho(\gamma)= \rho_{0}(\gamma)\circ
\phi$ for all $\gamma \in \Gamma$.  If there exists a group
element $\gamma \in \Gamma$ so that $\rho(\gamma)$ is partially
hyperbolic and so that the lift
$\widetilde{\mathcal{W}}_{\rho(\gamma)}^{u}$ of the unstable
foliation $\widetilde{\mathcal{W}}_{\rho(\gamma)}^{u}$ to
$\mathbb{R}^n$ is quasi-isometric, then $\phi$ is a homeomorphism.
\end{Prop}

\begin{Prf}
The map $\phi$ has degree one and is therefore surjective.  To
prove injectivity, it suffices to prove that $\phi$ is locally
injective.  It therefore suffices to show that a lift of $\phi$ to
$\mathbb{R}^n$ is locally injective.  First we argue that $\phi$
may be equivariantly lifted with respect to $\overline{\rho}$ and
$\pi$. Choose a lift $\tau:\mathbb{R}^n \rightarrow \mathbb{R}^n$
of $\phi$ and define $$\theta:\Gamma \times \mathbb{R}^n
\rightarrow \mathbb{R}^n$$ by
$\theta(\gamma,x)=\tau(\overline{\rho}(\gamma)(x))-\pi(\gamma)(\tau(x))$.
Since $\tau(\overline{\rho}(\gamma)(x))$ and
$\pi(\gamma)(\tau(x))$ project to the same point in the torus,
$\theta(\Gamma \times \mathbb{R}^n) \subset \mathbb{Z}^n$. As
$\mathbb{Z}^n$ is discrete, each $\theta(\gamma, \cdot)$ is
constant so that we may alternatively view $\theta$ as a map
$$\theta:\Gamma \rightarrow \mathbb{Z}^n \subset \mathbb{R}^n.$$
Viewed this way, $\theta$ is a one cocycle over $\pi$.  Indeed,
$$\theta(\gamma_1 \gamma_2)=\tau(\overline{\rho}(\gamma_1 \gamma_2)(0))
-\pi(\gamma_1 \gamma_2)(\tau(0))=$$
$$\tau(\overline{\rho}(\gamma_1)(\overline{\rho}(\gamma_2)(0)))
-\pi(\gamma_1)(\pi(\gamma_2)(\tau(0)))=$$
$$\theta(\gamma_1,\overline{\rho}(\gamma_2)(0)) +
\pi(\gamma_1)(\tau(\overline{\rho}(\gamma_2)(0))) -
\pi(\gamma_1)(\pi(\gamma_2)(\tau(0)))=$$
$$\theta(\gamma_1)+\pi(\gamma_1)(\theta(\gamma_2)),$$ for each $\gamma_1,\gamma_2 \in \Gamma.$
Since $H^1(\Gamma,\pi)=0$ (\cite[Theorem 2]{S}), there is a $v\in
\mathbb{R}^n$ so that
$$\theta(\gamma)=\pi(\gamma)(v)-v,$$ for each $\gamma \in \Gamma$.
 Define $$\overline{\phi}:\mathbb{R}^n \rightarrow \mathbb{R}^n$$
 by $\overline{\phi}(x) = \tau(x)+v$.  It is straightforward to
 check that $\overline{\phi}$ is equivariant with respect to
 $\overline{\rho}$ and $\pi$.  It therefore remains to show
 $\overline{\phi}$ is a cover of $\phi$.  First note that
 $\overline{\phi}$ descends to a map $\phi'$ homotopic to the identity
 on $\mathbb{T}^n$ since
 $$\overline{\phi}(x+z)=\tau(x+z)+v=\tau(x)+z+v=\overline{\phi}(x)+z,$$
 for each $x\in \mathbb{R}^n$ and $z\in \mathbb{Z}^n$.  Moreover, $\phi'$ is
  equivariant with respect to $\rho$ and $\rho_0$ and therefore coincides with $\phi$ by
 uniqueness.

 To
finish the argument, we show that $\overline{\phi}$ is locally
injective. Since $\phi$ is homotopic to the identity there is some
$M>0$ so that $\| \overline{\phi}(x)-x\|<M$.  If $\overline{\phi}$
is not locally injective, we may choose $x,y\in \mathbb{R}^n$ with
the same image and sufficiently close so that there exists a
piecewise $C^{1}$ curve $\sigma :[0,1] \rightarrow \mathbb{R}^n$
satisfying $\sigma(0)=x$, $\sigma(1) \in
\widetilde{\mathcal{W}}_{\rho(\gamma)}^{u}(y)$, and
$\dot{\sigma}\in E_{\rho(\gamma)}^{s}\oplus E_{\rho(\gamma)}^{c}$.
By equivariance,
$\overline{\phi}(\overline{\rho}(\gamma^n)(x))=\overline{\phi}(\overline{\rho}(\gamma^n)(y))$
for each natural number $n$.  Therefore,
$$2M> \|
\overline{\rho}(\gamma^{n})x-\overline{\rho}(\gamma^{n})y\| \geq
\|\overline{\rho}(\gamma^{n})y-\overline{\rho}(\gamma^{n})\sigma(1)\|
-\|\overline{\rho}(\gamma^{n})\sigma(1)-\overline{\rho}(\gamma^{n})x\|$$
$$\geq
\frac{1}{a}d_{\widetilde{\mathcal{W}}_{\gamma}^{u}}
(\overline{\rho}(\gamma^{n})
\sigma(1),\overline{\rho}(\gamma^{n})y) - K
\length(\overline{\rho}(\gamma^{n}) \sigma) - \frac{b}{a},$$ a
contradiction since the last term of this inequality grows
unbounded with $n$. \hfill $\Box$

\end{Prf}

\end{document}